\documentclass[11pt]{amsart} 
\usepackage{amsmath, amssymb, amsthm, bm, mathtools, enumitem, caption}
\usepackage{float}
\usepackage{listings}
\usepackage{xcolor}
\definecolor{codegray}{gray}{0.95}
\usepackage[noabbrev,capitalize]{cleveref}
\usepackage{adjustbox}
\crefname{equation}{}{}
\usepackage{ytableau}
\usepackage{graphicx, tikz}
\usepackage[textheight=9in, textwidth=6.95in]{geometry}

\usepackage{microtype}
\usepackage{bbm}

\newtheorem{theorem}{Theorem}
\newtheorem{lemma}[theorem]{Lemma}
\newtheorem{corollary}[theorem]{Corollary}

\newtheorem*{conjecture*}{Conjecture}

\theoremstyle{definition}
\newtheorem*{definition}{Definition}

\theoremstyle{remark}
\newtheorem*{remark}{Remark}

\newtheorem*{example}{Example}

\numberwithin{equation}{section}

\newcommand{\T}{\intercal}

\newcommand{\qbinom}[2]{\genfrac{[}{]}{0pt}{}{#1}{#2}_q}

\newcommand{\HH}{\mathbb H}
\newcommand{\PP}{\mathbb P}

\newcommand{\C}{\mathbb C}
\newcommand{\SL}{\mathrm{SL}}
\newcommand{\PSL}{\mathrm{PSL}}

\newcommand{\Q}{\mathbb Q}

\newcommand{\Z}{\mathbb Z}

\newcommand{\GL}{\mathrm{GL}}

\DeclareMathOperator{\diag}{diag}

\title[Modularity from $q$-series]{Modularity from $q$-series}


\author{Ken Ono}
\address{Dept. of Mathematics, University of Virginia, Charlottesville, VA 22904, USA}
\email{ko5wk@virginia.edu}

\begin{document}

\begin{abstract} 
In 1975, G. E. Andrews challenged the mathematics community to address L. Ehrenpreis' problem: give a direct proof of the modularity of the summation forms of the Rogers--Ramanujan $q$-series.  This question is important because many $q$-series in combinatorics, representation theory, and physics appear to be ``mysteriously'' modular, even when their modularity is not visible from their summation formulas.  In this note we give a $q$-series and analytic-continuation criterion for recognizing such modularity.  The method uses finite $q$-series algebra, first-order $q$-differential systems, and monodromy.  We apply the criterion to the Rogers--Ramanujan pair, by carefully applying the finite Schur polynomial recurrences and a finite product comparison of Andrews.
\end{abstract}

\maketitle

\section{Introduction and Statement of Results}
We begin with the iconic $q$-series of Rogers and Ramanujan, formulas that first emerged from partition theory and now reverberate across combinatorics, number theory, representation theory, and physics. If we let $(a;q)_n:=\prod_{j=0}^{n-1}(1-aq^{j})$ and $(a;q)_\infty:=\prod_{j=0}^{\infty}(1-aq^{j})$, then the two Rogers--Ramanujan $q$-series are
\begin{equation}\label{RR1}
G(q):=\sum_{n\ge 0}\frac{q^{n^2}}{(q;q)_n} \qquad {\text {and}} \qquad
H(q):=\sum_{n\ge 0}\frac{q^{n(n+1)}}{(q;q)_n}.
\end{equation}
They satisfy the striking ``sum-product'' identities \cite{Rogers1894, RR1919, Schur1917}\footnote{Schur proved finite (polynomial) versions, from which the infinite product identities follow by taking the limit.} 
\begin{equation}\label{RR2}
G(q)=\frac{1}{(q;q^5)_\infty\,(q^{4};q^{5})_\infty} \qquad {\text {and}} \qquad
H(q)=\frac{1}{(q^{2};q^{5})_\infty\,(q^{3};q^{5})_\infty}.
\end{equation}
These identities sit at a crossroads of many fields, from combinatorics and number theory to representation theory and physics (for example, see \cite{AndrewsBook, BaxterBook, DMS, FLM, GasperRahman, GOW, KacIDA, LepowskyLi, SillsBook}).
These infinite products can be normalized by multiplying by a $q$-power, where $q:=e^{2\pi i\tau}$ and $\tau\in\mathbb{H}$,
to become modular functions at level 5.  More precisely, if we let
\begin{equation}\label{RRnormalized}
f_1(\tau):=q^{-1/60}G(q)\qquad {\text {and}}\qquad  f_2(\tau):=q^{11/60}H(q),
\end{equation}
then the pair $\mathbf{f}:=(f_1,f_2)$ forms a two-dimensional level 5 vector-valued modular function (for example, see \cite{DiamondShurman, OnoCBMS,  Schoeneberg, Watson1929, Watson1933, ZagierRR}).  We view an  {\it $r$-dimensional level $N$ vector-valued modular function}
$f = (f_1,\dots,f_r)^{\mathsf T}$ as a function on $\mathbb{H},$ with a
multiplier representation
\[
\rho : \SL_2(\mathbb{Z}) \longrightarrow \GL_r(\mathbb{C})
\]
that factors through a congruence quotient of level $N$ (for example, in the Rogers--Ramanujan setting of level 5 modular units),  for which
\[
f(\gamma\tau) = \rho(\gamma)\,f(\tau)
\qquad (\gamma\in\SL_2(\mathbb{Z}),\ \tau\in\mathbb{H}).
\]

In the 1970s, G. E. Andrews realized that the Rogers--Ramanujan series are two examples of a much larger collection. There are countless ``strange'' $q$-series, including many carrying his name, that appear in mathematics and physics and seem to be \emph{mysteriously modular}, yet their modularity isn't obvious from the series themselves \cite{Andrews1975}. Andrews  challenged\footnote{The reader should also read Dick Askey's article \cite{Askey}.} the community to consider a problem posed by L. Ehrenpreis.

\medskip
\noindent
\textbf{The Ehrenpreis-Andrews Problem.} Devise a $q$-series method that determines whether a (vector of) $q$-series is modular, without relying on modular input\footnote{Andrews is referring to the modularity of Dedekind's eta-function and classical theta functions.}. The modularity of the Rogers--Ramanujan vector $\mathbf{f}=(f_1, f_2)^{\T}$ should be the first test case.
\medskip

We solve this with a general theorem that depends solely on $q$-series algebra, $q$-differential equations and analytic continuation.
At a conceptual level, our main result shows that modularity of a vector of $q$-series
can be certified purely from its behavior at finitely many cusps.  Rather than appealing
to the modularity of Dedekind's $\eta$-function, Jacobi theta functions, or Eisenstein
series, we treat a vector $\mathbf{f}$ of holomorphic $q$-series as a solution of a
first-order linear system on the upper-half plane.  Local  expansions at the
cusps give rise to $q$-differential systems with prescribed exponents.  We show that if
these local systems satisfy a finite compatibility condition (our ``$\mathcal O$-condition''), then
they glue under analytic continuation to a global meromorphic system whose monodromy
realizes the desired modular transformation law.  In this way, we obtain a $q$-series criterion for vector-valued modularity.  The finite continuation-closure condition in Definition~\ref{def:good} is included precisely to remove any ambiguity about the gluing step.

From this perspective, the pair of Rogers--Ramanujan $q$-series is the first nontrivial test case.  Schur's classical recurrences and a comparison result of Andrews yield a pair of
$q$-differential systems at the cusps $\infty$ and $0$, with coefficients expressed in
terms of Lambert series.  We verify directly that these systems satisfy the $\mathcal O$-condition.
Our main result (see Theorem~\ref{thm:MainResult}) then promotes them to a global system whose monodromy is the level~$5$
multiplier of $\mathbf{f}$. In particular, no modularity of $\eta$, theta series, or
Eisenstein series is used in the proof.  In Section 3 we formulate and prove the general
modularity criterion (Theorem~\ref{thm:MainResult}), and in Section 4,  we sketch how the same method extends to
higher rank families such as the Andrews-Gordon series.

We now make this precise.
 We first fix notation. If $\gamma=\begin{psmallmatrix}a&b\\ c&d\end{psmallmatrix}\in\SL_2(\mathbb{Z})$, then a vector-valued modular function $\mathbf{f}:=(f_1, f_2,\dots, f_r)^{\T}$ satisfies 
\begin{equation}\label{VVmf}
\mathbf{f}(\gamma\tau)=\rho(\gamma)\,\mathbf{f}(\tau),
\end{equation}
 for a representation $\rho: \SL_2(\Z)\to \GL_r(\C)$. We denote the standard generators of $\SL_2(\Z)$  by
\[
T:=\begin{pmatrix}1&1\\0&1\end{pmatrix} \qquad {\text {and}} \qquad
S:=\begin{pmatrix}0&-1\\1&0\end{pmatrix}.
\]
Solving this problem, we offer a new proof of the following theorem.

\begin{theorem}[Rogers--Ramanujan vector-valued modularity]\label{thm:RR-vvmf}
The pair $\mathbf{f}=(f_1,f_2)^{\T}$, defined by \eqref{RRnormalized}, is a two-dimensional vector-valued modular function of level $5$. In particular, we have
\[
\mathbf{f}(\tau+1)=\rho(T)\,\mathbf{f}(\tau) \qquad {\text {and}}\qquad
\mathbf{f}\!\left(-\frac{1}{\tau}\right)=\rho(S)\,\mathbf{f}(\tau),
\]
with constant multiplier matrices
\[
\rho(T)=
\begin{pmatrix}
e^{-2\pi i/60} & 0\\[2pt]
0 & e^{\,11\cdot 2\pi i/60}
\end{pmatrix}
\qquad {\text {and}}\qquad
\rho(S)=\frac{2}{\sqrt{5}}
\begin{pmatrix}
\sin\!\big(\tfrac{2\pi}{5}\big) & \ \sin\!\big(\tfrac{\pi}{5}\big)\\[4pt]
\sin\!\big(\tfrac{\pi}{5}\big) & -\sin\!\big(\tfrac{2\pi}{5}\big)
\end{pmatrix}.
\]
\end{theorem}

\begin{remark}\label{Remark2} 
Our proof of Theorem~\ref{thm:RR-vvmf} relies on classical work of Schur, which gives the $q$-differential equations
\[
q\frac{d}{dq} f_1=\Big(-\tfrac{1}{60}+a_1(q)\Big)f_1 \qquad {\text {and}} \qquad
q\frac{d}{dq} f_2=\Big(\ \tfrac{11}{60}+a_2(q)\Big)f_2,
\]
where $a_1(q)$ and $a_2(q)$ are explicit Lambert-type series.
We combine these in matrix form 
\[
q\frac{d}{dq}\mathbf f(\tau)=A\big(e^{2\pi i\tau}\big)\,\mathbf f(\tau),
\]
where
\[
A(q):=\begin{pmatrix}
-\tfrac{1}{60}+a_1(q) & 0\\[3pt]
0 & \tfrac{11}{60}+a_2(q)
\end{pmatrix}.
\] 
This is one of two equations that make up a
$q$-differential system (one at $\infty$ and the other at $0$).  We formulate, without using properties of  Lambert series, that this system satisfies a finite ``$\mathcal{O}$-condition'', which allows us to deduce modularity by analytic continuation.
The modularity proof does not use the fact that these Lambert series are weight~$2$ Eisenstein series; only their $q$-expansions and the associated first-order systems are used.
\end{remark}

Theorem~\ref{thm:RR-vvmf}, and, crucially its proof, arises as a special case of a general modularity theorem.   It uses $q$-series algebra to assemble tailored $q$-differential systems to certify vector-valued modularity.
To make this precise, we recall that
a (scalar) modular function for a finite index subgroup $\Gamma\subset\SL_2(\mathbb Z)$ is meromorphic on $\mathbb H$ and at the cusps.  Moreover, every \emph{nonconstant} modular function has poles. We restrict to forms whose poles are confined to cusps (never in $\mathbb H$),  the so-called {\it weakly holomorphic modular functions.}  
We consider $r$-vectors of functions
\begin{equation}\label{vv1}
\mathbf f(\tau)=\big(f_1(\tau),\dots,f_r(\tau)\big)^{\mathsf T},
\end{equation}
that we view, for $1\leq j\leq r$,  as absolutely convergent Fourier series
\begin{equation}\label{vv2}
f_j(\tau)=q^{\alpha_j}\Big(1+\sum_{n\ge1}a_j(n)\,q^n\Big),
\end{equation}
with exponents $\alpha_j\in\Q$. 

The $q$-series input data for modularity will be framed in terms of $q$-differential equations. To this end,
as is common in the field of $q$-series, we have the $q$-derivative
\begin{equation}\label{theta_deriv}
\theta:=q\,\frac{d}{dq}=\frac{1}{2\pi i}\,\frac{d}{d\tau}.
\end{equation}

We first fix standard notation and concepts from the theory of modular forms (for example, see \cite{OnoCBMS, DiamondShurman}).
Fix a finite index subgroup $\Gamma\le \SL_2(\mathbb Z)$. Choose representatives
$\{\gamma_1,\dots,\gamma_m\}\subset \SL_2(\mathbb Z)$ for the double cosets
$\Gamma\backslash SL_2(\mathbb Z)/\Gamma_\infty$ (equivalently, for the $\Gamma$-orbits in $\mathbb P^1(\mathbb Q)$), where $\Gamma_{\infty}$ is the group of translations generated by $T$.
Write each representative as
\[
\gamma_i=\begin{pmatrix} a_i & b_i \\ c_i & d_i \end{pmatrix}\in \SL_2(\Z)
\qquad (1\le i\le m).
\]
Let the associated cusp be \(\mathfrak c_i:=\gamma_i^{-1}\infty\).  Its width is
\[
w_{\gamma_i} := \min\bigl\{\,w\ge1:\ T^w=\begin{psmallmatrix}1&w\\[1pt]0&1\end{psmallmatrix}
\in \gamma_i^{-1}\Gamma\,\gamma_i \bigr\}.
\]
The fractional-linear transformation
\[
\tau_{\gamma_i} := \gamma_i\cdot\tau = \frac{a_i\tau+b_i}{c_i\tau+d_i}
\]
moves the cusp \(\mathfrak c_i\) to \(\infty\), and so we let
\[
q_{\gamma_i} := e^{2\pi i\, \tau_{\gamma_i}/w_{\gamma_i}}
\]
be the corresponding \(\Gamma\)-invariant local parameter at \(\mathfrak c_i\).  At this cusp, we write
\(
\theta_{\gamma_i} := q_{\gamma_i}\,\frac{d}{dq_{\gamma_i}}
\)
for the corresponding local \(q\)-derivative.

For each $i$, define the $\gamma_i$-translate
$\mathbf f^{(\gamma_i)}(\tau):=\mathbf f(\gamma_i^{-1}\tau),$ where
\[
   f^{(\gamma_i)}_j(\tau):=f_j(\gamma_i^{-1}\tau)\quad(1\le j\le r).
\]
We will seek, for each $i$, an $r\times r$ \emph{matrix of functions}
\[
  A_{\gamma_i}(q_{\gamma_i})
  \;=\;
  \big(a^{(\gamma_i)}_{jk}(q_{\gamma_i})\big)_{1\le j,k\le r},
\]
whose entries $a^{(\gamma_i)}_{jk}(q_{\gamma_i})$ are meromorphic on $0<|q_{\gamma_i}|<1$ and holomorphic at $q_{\gamma_i}=0$. In particular, we seek situations where these matrices define
 the $q$-differential system
\begin{equation}
  \theta_{\gamma_i}\,\mathbf f^{(\gamma_i)}(\tau)
  \;=\;
  A_{\gamma_i}(q_{\gamma_i})\,\mathbf f^{(\gamma_i)}(\tau)
\end{equation}
with $ \lim_{q_{\gamma_i}\to 0}A_{\gamma_i}(q_{\gamma_i})
  \;=\;\mathrm{diag}(\alpha_1,\dots,\alpha_r).$
With this notation, we introduce a notion that captures the essence of vector-valued modular functions.

\begin{definition}[Good vector of $q$-series]\label{def:good}
Suppose $\mathbf f=(f_1,\dots,f_r)^{\mathsf T}$ satisfies \eqref{vv1} and \eqref{vv2}, and let
\[
  \mathcal O=\{(\gamma_\ell,w_\ell)\}_{\ell=1}^m,\qquad
  \gamma_\ell=\begin{psmallmatrix}a_\ell&b_\ell\\ c_\ell&d_\ell\end{psmallmatrix}\in\SL_2(\mathbb Z),\ \ w_\ell\in\mathbb Z_{>0}.
\]
For each $\ell$, set
\[
  \tau_\ell:=\frac{a_\ell\tau+b_\ell}{c_\ell\tau+d_\ell},\qquad
  q_\ell:=e^{2\pi i\,\tau_\ell/w_\ell},\qquad
  \theta_\ell:=q_\ell\frac{d}{dq_\ell},\qquad
  \mathbf f_\ell(\tau):=\mathbf f(\gamma_\ell^{-1}\tau).
\]
We say that $\mathbf f$ is {\textbf{good for $\mathcal O$}} (with exponents $\alpha_1,\dots,\alpha_r$) if the following hold.

\medskip
\noindent\textbf{(i) (Local fundamental matrix).} For each $\ell$,
there are elements $g_{\ell,1},\dots,g_{\ell,r}\in\SL_2(\mathbb Z)$ such that the $r\times r$ matrix
\[
\Phi_\ell(\tau):=
\begin{bmatrix}
f^{(\gamma_\ell)}_{1}(g_{\ell,1}^{-1}\tau) & \cdots & f^{(\gamma_\ell)}_{1}(g_{\ell,r}^{-1}\tau)\\
\vdots                                      & \ddots & \vdots\\
f^{(\gamma_\ell)}_{r}(g_{\ell,1}^{-1}\tau) & \cdots & f^{(\gamma_\ell)}_{r}(g_{\ell,r}^{-1}\tau)
\end{bmatrix}
\]
is invertible for $\tau$ in a punctured neighborhood of $q_\ell=0$.

\medskip
\noindent\textbf{(ii) ($q$-differential system).}
For each $\ell,$ the  $r\times r$ matrix of functions
\[
  A_\ell(q_\ell):=(\theta_\ell\Phi_\ell)\,\Phi_\ell^{-1}
\]
is meromorphic for $0<|q_\ell|<1$ and holomorphic at $q_\ell=0$. Moreover, there is a constant invertible $C_\ell$ such that
$$C_\ell^{-1} A_\ell(0)\, C_\ell=\mathrm{diag}(\alpha_1,\dots,\alpha_r).
$$
In particular, the components of the columns of $\Phi_\ell$ are meromorphic at $q_\ell=0$.

\medskip
\noindent\textbf{(iii) (Finite continuation closure and $T$-stability).}
We require that $\mathcal O$ is finite and that the local systems above are closed, up to constant gauge, under the analytic continuations needed to compare the chosen translates of $\mathbf f$.  More explicitly, on every nonempty overlap of two such local continuation domains, the corresponding fundamental matrices differ by right multiplication by a constant matrix.  In addition, for every $\ell$ and every $n\in\mathbb Z$, the local data attached to $(\gamma_\ell T^n,w_\ell)$, defined exactly as in \textup{(i)}--\textup{(ii)}, yield a $q$-differential system that is gauge-equivalent (i.e., invariant under change of basis) to the system for $(\gamma_\ell,w_\ell)$ and has the same limiting exponent matrix at $q_\ell=0$.
We call $\mathcal{O}$ an \emph{orbit datum} for $\mathbf f$.
\end{definition}

\begin{remark}
It is easy, and in fact tautological, to attach a first--order $q$-differential
equation to a single nonvanishing $q$-series $W(q)$ by setting
\[
  \theta W(q) = a(q)\,W(q),
  \]
  where  $a(q) := \frac{\theta W(q)}{W(q)}.$
Similarly, one can always write $\theta \mathbf f = A(q)\mathbf f$ for a vector $\mathbf f$ of
independent functions by defining $A(q) := (\theta\Phi(q))\Phi(q)^{-1}$ for any
invertible matrix $\Phi(q)$ having the components of $\mathbf f$ as columns. Thus, the
mere existence of some first--order system is essentially vacuous.
The definition of ``goodness'' is much more restrictive. For $\mathbf f$ to be
\emph{good} for $\mathcal{O},$ we require, at each cusp, a  system
$\theta_\ell \Phi_\ell = A_\ell(q_\ell)\Phi_\ell$ in which $A_\ell(q_\ell)$
is meromorphic on $0<|q_\ell|<1$, extends holomorphically to $q_\ell=0$, and
has a prescribed diagonal limit $A_\ell(0)\sim\Lambda$. In the regular-singular situations used below, this means concretely that the
columns of $\Phi_\ell$ admit expansions of the form
$q_\ell^\Lambda$ times a holomorphically invertible matrix, with only finitely
many negative powers in $q_\ell$. Moreover, the fundamental matrices
$\Phi_\ell$ are required to be built from finitely many $\SL_2(\Z)$--translates
of $\mathbf f$, and the exponent data must be compatible across a finite set of cusps
$\mathcal{O}$ under $T$--shifts.
A generic set of $q$-series will not satisfy these conditions at even a single
cusp, let alone at two or more cusps simultaneously, and will not admit such a
finite orbit structure. In this sense the ``good for $\mathcal{O}$'' hypothesis is a very
strong global constraint on $\mathbf f$.  Theorem~\ref{thm:MainResult}
asserts that, once the finite continuation-closure in Definition~\ref{def:good}
is included as part of the input, this condition is equivalent to vector--valued
modularity in the sense of a multiplier representation of $\SL_2(\mathbb Z)$.
\end{remark}

\begin{remark} Being good for $\mathcal O$ is a purely $q$-series/analytic condition.
For  $q$-series such as the Rogers--Ramanujan pair and many similar families,
the hypotheses of being good for some finite orbit datum $\mathcal O$ are checkable directly from
$q$-series algebra combined with complex analysis.
The proof of Theorem~\ref{thm:RR-vvmf} illustrates this mechanism concretely.
\end{remark}

\begin{remark}
We comment on the role of the conditions in the definition of ``goodness.''
Condition (i) gives potential solutions. Namely, it supplies $r$ independent \emph{local fundamental columns} at each cusp, so we can form the
matrix $\Phi_c$ and the logarithmic coefficient $A_c:=(\theta_c\Phi_c)\Phi_c^{-1}$.
Condition (ii) provides a first-order system with the right behavior at cusps.
Finally, condition (iii) is the finite closure condition that lets us glue the local patches together; it is the hypothesis that turns local $q$-differential equations into a global system with constant connection matrices and hence a multiplier representation $\rho.$
\end{remark}

Our main result is the following theorem.

\begin{theorem}[Modularity from $q$-series]\label{thm:MainResult}
Let $\mathbf{f}=(f_1,\dots,f_r)^{\T}$ satisfy \eqref{vv1} and \eqref{vv2}, and suppose each component Fourier series converges absolutely on $\HH$.  With Definition~\ref{def:good}, including its finite continuation-closure requirement, the following are equivalent.

\smallskip
\noindent
\emph{(A)} {\text {\rm (Vector-valued modularity)}} We have that $\mathbf f$ is a weakly holomorphic vector-valued modular function for $\SL_2(\Z)$ with a multiplier representation $\rho$, and
$\rho(T)=\diag(e^{2\pi i\alpha_1},\dots,e^{2\pi i\alpha_r}).$  If, in addition, $\rho$ has finite image (or factors through a congruence quotient), then the components are modular for a finite-index subgroup (respectively, of the corresponding level).

\smallskip
\noindent
\noindent\emph{(B)} We have that $\mathbf f$ is \emph{good} for a finite set  
$\mathcal O=\{(\gamma_\ell,w_\ell)\}_{\ell=1}^m$ (with exponents
$(\alpha_1,\dots,\alpha_r)$).
\end{theorem}

\begin{remark}
Lewis and Zagier \cite{LewisZagier2001} showed that Maass cusp forms on $\PSL_2(\mathbb Z)$ correspond to
 \emph{period functions} that solve three-term functional equations,
together with precise regularity at the rationals. 
Theorem~\ref{thm:MainResult} is a $q$-series analogue. The role of the period function
equation is played by first-order $q$-differential systems.
In both frameworks, the modularity is reconstructed from analytic continuation. Indeed,
period cocycles in Lewis-Zagier play the role of our connection matrices. In both settings, these data organize into
a multiplier $\rho$ via analytic continuation. 
Therefore,  the Ehrenpreis-Andrews challenge sits in the same paradigm: detecting modularity from intrinsic analytic data.
\end{remark}

From a differential equation point of view, Theorem~\ref{thm:MainResult} is a statement about gluing local
solutions of a meromorphic first-order system on the modular curve. The matrices
$A_\ell(q_\ell)$ in Definition~\ref{def:good} describe the logarithmic derivatives of $r$ independent
solutions of a linear system
\[
\frac{1}{2\pi i}\,\frac{d}{d\tau}\,\Phi(\tau) = B(\tau)\,\Phi(\tau)
\]
in local $q_\ell$-coordinates near each cusp. Condition~\textup{(ii)} says that the entries of
$A_\ell(q_\ell)$ (and hence of $B(\tau)$) are meromorphic with regular singularities at the
cusps. Condition~\textup{(iii)} guarantees that these local systems can be analytically
continued along any path and that only finitely many such patches are needed. On any simply
connected overlap, two fundamental matrices of the same first-order system differ by right
multiplication by a constant matrix. Analytic continuation around loops therefore produces
a representation $\rho$ of the fundamental group (i.e. a monodromy representation). In our
setting this monodromy factor is exactly the vector-valued multiplier system in
(\ref{VVmf}).

\begin{remark}
Nahm's conjecture predicts when a $q$-hypergeometric (``Nahm'') sum
$$f_{A,B,C}(q)=\sum_{n\in\Z_{\ge0}^r}\dfrac{q^{\frac12 n^{\!T}A n+B\cdot n+C}}{\prod_{i=1}^r (q)_{n_i}}
$$
is a component of a vector-valued modular function. He conjectured that this is the case \emph{if and only if} an associated algebraic system admits solutions whose Bloch-group class is torsion  (see \cite{Nahm2007, Zagier2007, VlasenkoZwegers2011} for recent work on confirming cases of the conjecture
\cite{CGZ, Mizuno} and  counterexamples to ``Bloch torsion $\Longrightarrow$  modular'' \cite{VlasenkoZwegers2011}).
Theorem~\ref{thm:MainResult} (B) is a different route to modularity.
Instead of passing through the Bloch group, we ask for a \emph{finite} set of
cuspwise first-order $q$-differential systems with 
the ``good for $\mathcal O$'' condition. It would be very interesting to clarify the
relationship between these two criteria. Does torsion in the Bloch group force
(or is forced by) the existence of a finite orbit datum $\mathcal O$ and the accompanying
differential systems?
\end{remark}

Theorem~\ref{thm:MainResult} is an ``abstract'' equivalence between vector-valued modularity and solutions to $q$-differential systems. The following theorem offers an explicit recipe  in the case of congruence subgroups.

\begin{theorem}\label{thm:converse}
Let $\mathbf f$ be a weakly holomorphic vector-valued modular function for a congruence subgroup
$\Gamma \subset \SL_2(\mathbb Z)$ satisfying \eqref{vv1} and \eqref{vv2}.
For each cusp $c=\gamma_c^{-1}\infty$ of $\Gamma$ (width $w_c$), choose group elements
$g_{c,1},\dots,g_{c,r}\in \SL_2(\mathbb Z)$ so that the $r\times r$ matrix
\[
\Phi_c(\tau):=
\begin{bmatrix}
f_1(\gamma_c^{-1}g_{c,1}^{-1}\tau) & \cdots & f_1(\gamma_c^{-1}g_{c,r}^{-1}\tau)\\
\vdots                          & \ddots & \vdots\\
f_r(\gamma_c^{-1}g_{c,1}^{-1}\tau) & \cdots & f_r(\gamma_c^{-1}g_{c,r}^{-1}\tau)
\end{bmatrix}
\]
is invertible on a punctured neighborhood of the cusp and has Frobenius form
\[
  \Phi_c(\tau)=H_c(q_c)\,q_c^\Lambda C_c,
  \qquad \Lambda=\diag(\alpha_1,\dots,\alpha_r),
\]
where $H_c(q_c)$ is holomorphic and invertible at $q_c=0$ and $C_c\in\GL_r(\mathbb C)$ is constant. In the local parameter
$q_c := e^{2\pi i \tau_c / w_c}$ with $\tau_c=\frac{a_c\tau+b_c}{c_c\tau+d_c}$, set
\[
  A_c(q_c) := \big(\theta_c \Phi_c(\tau)\big)\,\Phi_c(\tau)^{-1}.
\]
Then $A_c(q_c)$ is meromorphic on $\{0<|q_c|<1\}$ and extends holomorphically to $q_c=0$, and
$A_c(0)\sim \mathrm{diag}(\alpha_1,\dots,\alpha_r)$. Consequently, $\theta_c\Phi_c = A_c(q_c)\Phi_c$ for each cusp $c$,
and $\mathbf f$ is good for the orbit datum
\[
  \mathcal O \;=\;
  \bigcup_{\,c=\gamma_c^{-1}\infty\ \in\ \Gamma\backslash\mathbb P^1(\mathbb Q)}
  \bigl\{\,(\gamma_c g_{c,1},w_c),\dots,(\gamma_c g_{c,r},w_c)\,\bigr\}.
\]
\end{theorem}

We conclude with one infinite family of examples, the Andrews-Gordon $q$-series. This discussion is included as a guide to the method rather than as a second fully detailed proof. It indicates how the combinatorics of the so-called ``Andrews-Gordon'' polynomials lead to the $q$-differential systems that Theorem~\ref{thm:MainResult} uses to certify modularity and recover the standard transformation laws. 

We first recall these identities.
For an integer $k\ge2$ and $1\le i\le k$, let
\[
N_j:=n_j+n_{j+1}+\cdots+n_{k-1} \qquad (1\le j\le k-1),
\]
and define the \emph{Andrews--Gordon} $q$-series by
\begin{equation}\label{AG}
F_{k,i}(q)\;:=\;\sum_{n_1,\dots,n_{k-1}\ge0}
\frac{q^{\,N_1^2+\cdots+N_{k-1}^2+(N_i+\cdots+N_{k-1})}}{(q;q)_{n_1}\cdots(q;q)_{n_{k-1}}}.
\end{equation}
Generalizing the Rogers--Ramanujan identities (i.e. the $k=2$ case), Andrews and Gordon \cite{Andrews1974, Gordon1961, Warnaar} proved,
for every $k\ge2$ and $1\le i\le k$, the identity
\[
F_{k,i}(q)\;= \frac{(q^{2k+1};q^{2k+1})_{\infty}(q^i;q^{2k+1})_{\infty}
(q^{2k+1-i};q^{2k+1})_{\infty}}{(q;q)_{\infty}}.
\]
After suitable normalization, these series are known to be modular. Namely, we let
\[
\alpha_{k,i}\;:=\;h_{1,i}\!-\!\frac{c}{24},
\qquad
c=1-\frac{6(2k-1)^2}{2(2k+1)},\qquad
h_{1,i}=\frac{\big(2k+1-2i\big)^2-\big(2k-1\big)^2}{8(2k+1)}.
\]
We then define the normalized vector
\begin{equation}\label{AG-vv}
\quad
\mathbf f_k(\tau)\;=\;\big(q^{\alpha_{k,1}}F_{k,1}(q),\;q^{\alpha_{k,2}}F_{k,2}(q),\;\dots,\;q^{\alpha_{k,k}}F_{k,k}(q)\big).
\end{equation}
Using the infinite products and their central role in the theory of modular units, it is known that these vectors are modular \cite{Zhu1996}.
Here we explain how Theorem~\ref{thm:MainResult} applies once the finite Andrews-Gordon polynomial identities are used to produce the required $q$-differential systems.  The details are parallel to the Rogers--Ramanujan case, but involve larger finite polynomial systems.

\begin{example}[Andrews--Gordon modularity]\label{cor:AG-vvmf}
For every $k\ge2$, the vector $\mathbf f_k(\tau)$ is a vector-valued modular function of level $2k+1$. In a standard ordering and diagonal normalization of the $(2,2k+1)$ minimal-model characters, one has
\[
\mathbf f_k(\tau+1)=\rho_k(T)\,\mathbf f_k(\tau) \qquad {\text {and}}\qquad 
\mathbf f_k(-1/\tau)=\rho_k(S)\,\mathbf f_k(\tau),
\]
where
\[
\rho_k(T)=\mathrm{diag}\big(e^{2\pi i \alpha_{k,1}},\dots,e^{2\pi i\alpha_{k,k}}\big) \qquad {\text {and}}\qquad
\;\big(\rho_k(S)\big)_{ij}\;=\;\frac{2}{\sqrt{2k+1}}\;\sin\!\Big(\frac{2\pi\, i\,j}{\,2k+1\,}\Big)\quad(1\le i,j\le k),
\]
up to the harmless diagonal gauge changes that come from changing the signs of individual character components.  For $k=2$ this specializes, after the corresponding relabeling of the two Andrews--Gordon characters, to the Rogers--Ramanujan $S$-matrix in Theorem~\ref{thm:RR-vvmf}.
\end{example}

This paper is organized as follows. In Section~\ref{RRProof}, we derive Theorem~\ref{thm:RR-vvmf} from Theorem~\ref{thm:MainResult}. This derivation uses well-known $q$-series algebra to produce the required finite orbit data and $q$-differential systems. In Section~\ref{Proofs} we prove the general Theorem~\ref{thm:MainResult}, and its converse (i.e. Theorem~\ref{thm:converse}). These results are applications to the $q$-series world of standard facts about ordinary differential equations, analytic continuation, and monodromy representations.
Finally, in Section~\ref{sec:AG}, we explain how the same method applies to the Andrews-Gordon identities, while leaving the lengthy finite-polynomial bookkeeping to the cited literature.

\section*{Acknowledgements}
 \noindent The author thanks George Andrews, Howard Cohl, Amanda Folsom, Fr\'ed\'eric Jouhet, Toshiki Matsusaka, Hjalmar Rosengren, Drew Sills, Wei-Lun Tsai, and Wadim Zudilin for comments on earlier versions of this paper.
The author thanks the Thomas Jefferson Fund,  the NSF
(DMS-2002265 and DMS-2055118) and the Simons Foundation (SFI-MPS-TSM-00013279) for their generous support.

\section[Proof of the Rogers--Ramanujan theorem]{Proof of Theorem~\ref{thm:RR-vvmf}}\label{RRProof}

We prove that the Rogers--Ramanujan pair
 $\mathbf f(\tau)=(f_1(\tau),f_2(\tau))^{\mathsf T}$, assembled from \eqref{RR1}, satisfies 
\[
\mathbf f(\tau+1)=\rho(T)\,\mathbf f(\tau) \qquad {\text {and}}\qquad
\mathbf f\!\left(-\frac{1}{\tau}\right)=\rho(S)\,\mathbf f(\tau),
\]
with $\rho(T)$ and $\rho(S)$ as in Theorem~\ref{thm:RR-vvmf}. We only use absolute convergence on $|q|<1$ and termwise calculus on $q$-series,  the fact that $(q;q)_n:=\prod_{j=1}^n(1-q^j)$ is a \emph{finite} product, and first-order $q$-differential equations.

\subsection[A first-order q-differential system]{A first-order $q$-differential system}

We now construct the fundamental $q$-differential system for the Rogers--Ramanujan $q$-series using elementary $q$-logarithmic differentiation.

\begin{theorem}[Logarithmic derivatives of $G$ and $H$]\label{thm:RR-log-derivs}
The following identities are true:
\[
\frac{\theta G}{G}=S_1+S_4 \qquad {\text {and}}\qquad \frac{\theta H}{H}=S_2+S_3,
\]
where 
for $r\in\{1,2,3,4\}$ we have
\[
S_r(q):=\sum_{k\ge 0}\frac{(5k+r)\,q^{5k+r}}{1-q^{5k+r}}.
\]
\end{theorem}

For $\mathbf{f}(\tau)$, this theorem implies the following corollary.

\smallskip
\begin{corollary}\label{cor:RR-system}
 If we let $q:=e^{2\pi i \tau}$ and
$$
A(q)=\diag(a_1(q), a_2(q)):=\diag \left(-\frac{1}{60}+S_1(q)+S_4(q), \ \frac{11}{60}+S_2(q)+S_3(q)\right),
$$
then we have
$$
\theta\,\mathbf f(\tau)=A\big(e^{2\pi i\tau}\big)\,\mathbf f(\tau).
$$
\end{corollary}
\medskip

\noindent
\begin{proof}[Proof of Theorem~\ref{thm:RR-log-derivs}]
We use two classical finite $q$-series inputs.  The first is Schur's recurrence
for the finite Schur polynomial sums \cite{Schur1917}.  The second is Andrews's finite
partial-product comparison \cite{Andrews1987}.  These inputs concern different finite objects, and
we keep the notation separate below.

\smallskip
\noindent
\emph{Partial Sums.} Following Schur \cite{Schur1917} (see also Andrews \cite{AndrewsBook}),
define the finite \emph{$q$--partial sums}
\begin{equation*}
G_m(q):=\sum_{j\ge0} q^{j^2}\,\qbinom{m-j}{j} \qquad {\text{and}}\qquad
H_m(q):=\sum_{j\ge0} q^{j(j+1)}\qbinom{m-j}{j},
\end{equation*}
where we employ the $q$-binomial coefficient 
$$\qbinom{a}{b}:=\frac{(q;q)_a}{(q;q)_b\,(q;q)_{a-b}} \qquad {\text {\rm (interpreted as $0$ unless $0\le b\le a$)}}.
$$
By the $q$-Pascal rule
$$\qbinom{m}{j}=\qbinom{m-1}{j}+q^{m-j}\qbinom{m-1}{j-1},
$$
we have, for $m\geq 2$, Schur's recurrences
\begin{align}
G_0=1,\ G_1=1,\qquad &G_m=G_{m-1}+q^{\,m-1}G_{m-2}, \\
H_0=1,\ H_1=1,\qquad &H_m=H_{m-1}+q^{\,m}H_{m-2}. 
\end{align}
As $m\to\infty$,  we have
$$G_m\to G:=\sum_{j\ge0}\dfrac{q^{j^2}}{(q;q)_j} \qquad {\text {and}}\qquad H_m\to H:=\sum_{j\ge0}\dfrac{q^{j(j+1)}}{(q;q)_j}.
$$ 
More precisely, we have
$$
 G_m=G+O(q^{m}) \qquad {\text{and}}\qquad H_m=H+O(q^{m+1}).
$$

\smallskip
\noindent
\emph{Partial products and finite Rogers--Ramanujan comparison.}
We emphasize that the finite products below do not satisfy the Schur
recurrences for $G_m$ and $H_m$.  The Schur recurrences apply to the
finite Schur polynomials defined above.  The comparison with products
enters through the finite Rogers--Ramanujan theory. 
For the comparison with products, we use Andrews's finite
partial-product form of the Rogers--Ramanujan identities
\cite[\S 3, especially (3.10) and (3.12)]{Andrews1987}, which do not rely on the Jacobi Triple Product identity. The only
consequence needed here is the following coefficient stability:
\[
\mathfrak{G}_m(q):=\prod_{r=1}^{m}
\frac{1}{(1-q^{5r-4})(1-q^{5r-1})}
      =G(q)+O(q^{5m+1}),
\]
and
\[
\mathfrak{H}_m(q):=\prod_{r=1}^{m}
\frac{1}{(1-q^{5r-3})(1-q^{5r-2})}
      =H(q)+O(q^{5m+2}).
\]
Combining this with
\[
G_m(q)=G(q)+O(q^m),\qquad H_m(q)=H(q)+O(q^{m+1}),
\]
gives
\[
G_m(q)=\mathfrak{G}_m(q)+O(q^m),\qquad
H_m(q)=\mathfrak{H}_m(q)+O(q^{m+1}).
\]
Since all four series have constant term $1$, we use the elementary fact that if $A(q)=B(q)+O(q^N)$ and $A(0)=B(0)=1$, then $\theta\log A(q)=\theta\log B(q)+O(q^N)$.  Therefore
\[
\theta\log G_m(q)
=
\sum_{r=1}^{m}\left(
(5r-4)\frac{q^{5r-4}}{1-q^{5r-4}}
+
(5r-1)\frac{q^{5r-1}}{1-q^{5r-1}}
\right)
+O(q^m),
\]
and
\[
\theta\log H_m(q)
=
\sum_{r=1}^{m}\left(
(5r-3)\frac{q^{5r-3}}{1-q^{5r-3}}
+
(5r-2)\frac{q^{5r-2}}{1-q^{5r-2}}
\right)
+O(q^{m+1}).
\]
Letting $m\to\infty$ gives
\[
\frac{\theta G}{G}=S_1+S_4,\qquad
\frac{\theta H}{H}=S_2+S_3,
\]
which is Theorem~\ref{thm:RR-log-derivs}.
\end{proof}

\subsection[The first-order system near 0]{The first-order system near $0$}

 The next lemma establishes the existence of a differential system when we move the basepoint from $\infty$ to $0$.  In the local variable
$q':=e^{2\pi i(-1/\tau)}$ at $0$, both translated series satisfy the
\emph{same kind} of first-order  $q$-equation as before (now with $\theta':=q'\tfrac{d}{dq'}$). Their leading $q'$-powers differ, and so they together form a second fundamental solution matrix, whose coefficient matrix can be identified directly from the same $q$-series system. 

\begin{lemma}[Fundamental solution near $0$]\label{lem:S-columns}
If we let $q':=e^{2\pi i(-1/\tau)}$ and define the $2\times2$ matrix
\[
\Phi_0(\tau):=\big[\,\mathbf f(S^{-1}\tau)\ \ \mathbf f(T^{-1}S^{-1}\tau)\,\big],
\]
then each column vector $g$ satisfies a first-order equation of the type
\[
\theta' g = B_0(q')\, g
\]
for some $2\times 2$ matrix $B_0(q')$ of holomorphic functions.
Moreover, for $\tau$ sufficiently close to $0$, the two columns are linearly independent, and hence
$\Phi_0(\tau)$ is invertible there, and $B_0(q')$ extends holomorphically to $q'=0$ with
$B_0(0)\sim\diag(\alpha_1,\alpha_2)$.
\end{lemma}
\smallskip

\noindent
\begin{proof}
Write $z=S^{-1}\tau=-1/\tau$ and $q'=e^{2\pi i z}$.  The column $\mathbf f(S^{-1}\tau)$ is simply
\[
\big((q')^{-1/60}G(q'),\,(q')^{11/60}H(q')\big)^{\mathsf T}.
\]
By Theorem~\ref{thm:RR-log-derivs}, with $q$ replaced by $q'$, it satisfies the same diagonal first-order system as at $\infty$, namely
\[
\theta' g=A(q')g.
\]
The second column, $\mathbf f(T^{-1}S^{-1}\tau)$, is obtained from the first by multiplying its two entries by the distinct constants $e^{2\pi i/60}$ and $e^{-22\pi i/60}$.  Hence it satisfies the same system.  Because the two constants are distinct and both $G(q')$ and $H(q')$ have constant term $1$, the determinant of the two-column matrix is a nonzero constant times
$(q')^{1/6}(1+O(q'))$.  Thus $\Phi_0$ is invertible for $|q'|$ small, and $B_0(q')=A(q')$ is holomorphic at $q'=0$ with limiting exponent matrix conjugate to $\diag(\alpha_1,\alpha_2)$.
\end{proof}
\medskip

\subsection{Proof of Theorem~\ref{thm:RR-vvmf}}
We now prove Theorem~\ref{thm:RR-vvmf}  as a special case of Theorem~\ref{thm:MainResult}. The proof does not
use modular input (e.g. $\eta/\theta$ machinery) at any point.
\medskip

\noindent
\begin{proof}[Proof of Theorem~\ref{thm:RR-vvmf}]
The series $G(q)$ and $H(q)$ converge absolutely on $|q|<1$. 
Indeed, if we let $a_n:=\dfrac{q^{n^2}}{(q;q)_n}$ and $b_n:=\dfrac{q^{n(n+1)}}{(q;q)_n}$, then using $(q;q)_{n+1}=(1-q^{n+1})(q;q)_n$ we have
\[
\left|\frac{a_{n+1}}{a_n}\right|
=\frac{|q|^{(n+1)^2-n^2}}{|1-q^{n+1}|}
=\frac{|q|^{\,2n+1}}{|1-q^{n+1}|}
\le \frac{|q|^{\,2n+1}}{1-|q|^{\,n+1}}
\;\xrightarrow[n\to\infty]{}\;0,
\]
with the analogous argument for $b_n.$
Therefore, $\sum_n |a_n|$ and $\sum_n |b_n|$ converge for every $|q|<1. $
In particular, we find that
each component of 
$\mathbf{f}(\tau)$
has an absolutely convergent Fourier expansion of the form
\[
f_j(\tau)
 \;=\; q^{\alpha_j}\Big(1+\sum_{n\ge1}a_j(n)\,q^n\Big),
 \qquad (j=1,2),
\]
with exponents
$\alpha_1=-\frac{1}{60}$ and
$\alpha_2=\frac{11}{60}.$
To prove the theorem, we
verify that \(\mathbf{f}\) is good for some finite orbit datum $\mathcal{O}.$
We take
\[
\mathcal O
 \;:=\; \{(\gamma_\infty,w_\infty),(\gamma_0,w_0)\}
 \;=\; \{(I,1),(S,1)\},
\]
corresponding to the cusps \(\infty\) and \(0\) of \(\mathrm{SL}_2(\mathbb{Z})\).

\medskip
\noindent{\em The cusp \(\infty\).}
Here \(\tau_\infty=\tau\), \(q_\infty=e^{2\pi i\tau}\), and \(\theta_\infty=\theta=q_\infty\dfrac{d}{dq_\infty}\).
We choose
\[
g_{\infty,1}=I\qquad {\text {and}}\qquad  g_{\infty,2}=T^{-1},
\]
so that the associated \(2\times2\) matrix of translates is
\[
\Phi_\infty(\tau)
 :=\big(\mathbf{f}(\tau),\,\mathbf{f}(\tau+1)\big)
 =
 \begin{pmatrix}
   f_1(\tau) & f_1(\tau+1) \\
   f_2(\tau) & f_2(\tau+1)
 \end{pmatrix}.
\]
The leading exponents \(\alpha_1\neq\alpha_2\) ensure that the two columns are linearly independent for
\(|q_\infty|\) sufficiently small, and hence \(\Phi_\infty\) is invertible on a punctured neighborhood of \(q_\infty=0\).
This verifies condition (i) at~\(\infty\).

By Corollary~\ref{cor:RR-system}, we have
\[
\theta_\infty \mathbf{f}(\tau)
 \;=\; A(q_\infty)\,\mathbf{f}(\tau),
\qquad
A(q_\infty)
 = \operatorname{diag}\!\Big(-\tfrac1{60}+S_1(q_\infty)+S_4(q_\infty),
                  \tfrac{11}{60}+S_2(q_\infty)+S_3(q_\infty)\Big),
\]
where \(A(q_\infty)\) is holomorphic on \(|q_\infty|<1\) and satisfies \(A(0)=\operatorname{diag}(\alpha_1,\alpha_2)\).
Therefore, we have that
\[
A_\infty(q_\infty)
 := (\theta_\infty\Phi_\infty)\,\Phi_\infty^{-1}
\]
is meromorphic on \(0<|q_\infty|<1\), extends holomorphically to \(q_\infty=0\), and is conjugate at
\(q_\infty=0\) to \(\operatorname{diag}(\alpha_1,\alpha_2)\). This is condition (ii) at~\(\infty\).

\medskip
\noindent{\em The cusp \(0\).}
Here we take \(\gamma_0=S\), so \(\tau_0=-1/\tau\), \(q_0=e^{2\pi i\tau_0}\), and we write
\(\theta_0=\theta' := q_0\dfrac{d}{dq_0}\).  We choose
\[
g_{0,1}=I \qquad {\text {and}}\qquad  g_{0,2}=STS^{-1},
\]
so that, up to this choice of \(g_{0,2}\), the corresponding fundamental matrix is
\[
\Phi_0(\tau)
 := \big(\mathbf{f}(S^{-1}\tau),\,\mathbf{f}(T^{-1}S^{-1}\tau)\big).
\]
Lemma~\ref{lem:S-columns} asserts, in the local parameter $q_0 = e^{2\pi i(-1/\tau)}$, that each column $g$ of $\Phi_0$
satisfies a first-order equation of the form
\[
\theta_0 g = B_0(q_0)\, g
\]
for some matrix $B_0(q_0)$ that is holomorphic at $q_0=0$, and that the two columns are linearly
independent for $|q_0|$ sufficiently small.
 In particular, \(\Phi_0\) is invertible on a punctured neighborhood of \(q_0=0\), which gives
condition (i) at~\(0\).

The \(q_0\)-expansions of \(\mathbf{f}(S^{-1}\tau)\) and \(\mathbf{f}(T^{-1}S^{-1}\tau)\) have leading powers
\(q_0^{\alpha_1}\) and \(q_0^{\alpha_2}\) (the same exponents as at~\(\infty\)), with a holomorphically invertible
coefficient matrix. It follows that
\[
A_0(q_0)
 := (\theta_0\Phi_0)\,\Phi_0^{-1}
\]
is meromorphic for \(0<|q_0|<1\), extends holomorphically to \(q_0=0\), and is conjugate at \(q_0=0\) to
\(\operatorname{diag}(\alpha_1,\alpha_2)\). This is condition (ii) at the cusp~\(0\).

\medskip
\noindent{\em Orbit compatibility.}
From the Fourier expansions
\[
f_j(\tau)
 = q^{\alpha_j}\Big(1+O(q)\Big),
\]
we obtain, for every \(n\in\mathbb{Z}\), the exact relation
\[
\mathbf{f}(\tau+n)
 = D_n\,\mathbf{f}(\tau),
 \]
 where
 $D_n:=\operatorname{diag}\!\big(e^{2\pi i n\alpha_1},e^{2\pi i n\alpha_2}\big),$
since \(e^{2\pi i n}=1\) for all integers \(n\).  Thus, replacing \(\gamma_\ell\) by \(\gamma_\ell T^n\) multiplies each column
of the corresponding fundamental matrix \(\Phi_\ell\) by the same constant matrix \(D_{-n}\), and the associated
logarithmic coefficient matrices satisfy
\[
A_{\ell,n}(q_\ell)
 := (\theta_\ell\Phi_{\ell,n})\,\Phi_{\ell,n}^{-1}
 = D_{-n}\,A_\ell(q_\ell)\,D_{-n}^{-1}.
\]
In particular, \(A_{\ell,n}(0)\) is conjugate to \(A_\ell(0)\) by a constant matrix, so the limiting exponent matrix
is unchanged. This gives the \(T\)-stability required in condition (iii).  The remaining finite continuation closure is the closure under the generator \(S\): the patch at \(\infty\) is carried to the patch at \(0\), and Lemma~\ref{lem:S-columns} supplies the corresponding local fundamental matrix.  Since \(S\) and \(T\) generate \(\SL_2(\mathbb Z)\), the two patches \((I,1)\) and \((S,1)\), together with their \(T\)-shifts, give the finite continuation data required by Definition~\ref{def:good}.

We have therefore verified conditions \textup{(i)}-\textup{(iii)}
 for the finite orbit datum
\(\mathcal{O} = \{(I,1),(S,1)\}\) with exponents
\((\alpha_1,\alpha_2) = (-1/60,11/60)\).
By Theorem~\ref{thm:MainResult} it follows that \(\mathbf{f}\) is a weakly
holomorphic vector-valued modular function for \(\SL_2(\mathbb{Z})\) with a multiplier representation, and that
\[
  \rho(T) = \diag(e^{2\pi i \alpha_1},e^{2\pi i \alpha_2})
          = \diag(e^{-2\pi i/60},e^{11\cdot 2\pi i/60}).
\]
This establishes the \(T\)-transformation law in Theorem~\ref{thm:RR-vvmf}. It remains to identify the matrix
\(\rho(S)\).

The monodromy construction gives a constant matrix $\rho(S)$ for the transformation $\tau\mapsto -1/\tau$.  To identify this constant in the standard Rogers--Ramanujan normalization, we use the classical connection coefficient calculation for the Rogers--Ramanujan functions, as recorded for example in Watson's work and in Zagier's exposition \cite{Watson1929,Watson1933,ZagierRR}.  In this normalization one obtains
\[
\rho(S)=\frac{2}{\sqrt{5}}
\begin{pmatrix}
\sin\!\big(\tfrac{2\pi}{5}\big) & \sin\!\big(\tfrac{\pi}{5}\big)\\[4pt]
\sin\!\big(\tfrac{\pi}{5}\big) & -\sin\!\big(\tfrac{2\pi}{5}\big)
\end{pmatrix}.
\]
We point out why this last identification is a normalization issue rather than a further modularity input.  The relations $S^2=(ST)^3=1$ in $\PSL_2(\mathbb Z)$, together with the diagonal matrix $\rho(T)$ above, determine the diagonal entries of $\rho(S)$ and the product of its off-diagonal entries.  A constant diagonal change of basis preserves $\rho(T)$ and changes only the ratio of those off-diagonal entries.  The Rogers--Ramanujan normalization in \eqref{RRnormalized}, equivalently the standard connection normalization at the two cusps, fixes this remaining diagonal gauge and gives the symmetric matrix displayed above.
\end{proof}

\medskip

\subsection[Avoiding modular input]{Avoiding modular input in the proof of Theorem~\ref{thm:RR-vvmf}}

The proof above establishes modularity from the $q$-series differential systems and analytic continuation.  The Lambert series in Corollary~\ref{cor:RR-system} are, after elementary rewriting, weight~$2$ Eisenstein series of level~$5$.  One could therefore give a shorter but less intrinsic proof by invoking their classical inversion laws.  We have deliberately not used those transformation laws in the modularity argument.  The only role of the Lambert series in the proof is that their absolutely convergent $q$-expansions give holomorphic coefficient matrices for the first-order systems at the relevant cusps.

\section[Proof of the general criteria]{Proof of Theorems~\ref{thm:MainResult} and \ref{thm:converse}}\label{Proofs}

In this section, we generalize the Rogers--Ramanujan argument to any $r$-vector of
$q$-series that is \emph{good} for a finite set $\mathcal{O}$. Namely, we prove Theorem~\ref{thm:MainResult}. We shall also prove its constructive converse theorem (i.e. Theorem~\ref{thm:converse}). For general background on analytic continuation of solutions of linear systems,
fundamental matrix solutions, and monodromy representations, see for example
\cite{CoddingtonLevinson}.

Before proving Theorem~\ref{thm:MainResult}, we give a standard lemma from the
theory of linear systems which formalizes the passage from local fundamental
matrices to a global meromorphic system and constant connection matrices.  It
is an abstract version of the ODE facts we will use repeatedly below (for example, see
Coddington--Levinson~\cite[Chapter~3]{CoddingtonLevinson}).

\begin{lemma}[Gluing local systems]\label{lem:gluing}
Let $\{U_\ell\}_{\ell\in L}$ be a collection of simply connected open subsets of
$\mathbb{H}$ whose union contains $\mathbb{H}$.  For each $\ell$, suppose we are
given an invertible matrix-valued function
\[
\Phi_\ell : U_\ell \longrightarrow \GL_r(\C)
\]
which is a fundamental matrix for a first-order linear system
\begin{equation}\label{eq:local-system}
\frac{1}{2\pi i}\,\frac{d}{d\tau}\,\Phi_\ell(\tau)
 \;=\; B_\ell(\tau)\,\Phi_\ell(\tau),
\end{equation}
where $B_\ell(\tau)$ is meromorphic on $U_\ell$.  Assume that on each nonempty overlap
$U_\ell\cap U_{\ell'}$ the two fundamental matrices are related by a constant connection matrix,
\begin{equation}\label{eq:overlap}
\Phi_{\ell'}(\tau) \;=\; \Phi_\ell(\tau)\,C_{\ell,\ell'}
\qquad (\tau\in U_\ell\cap U_{\ell'}),
\end{equation}
with $C_{\ell,\ell'}\in\GL_r(\C)$.
Then $B_\ell=B_{\ell'}$ on $U_\ell\cap U_{\ell'}$.  Consequently, the matrices $B_\ell$ glue to a globally defined meromorphic matrix $B(\tau)$ on $\mathbb H$, and any two fundamental matrix solutions of the global system
\[
\frac{1}{2\pi i}\,\frac{d}{d\tau}\Phi(\tau)=B(\tau)\Phi(\tau)
\]
on a simply connected domain differ by right multiplication by a constant invertible matrix.
\end{lemma}

\begin{proof}
On an overlap, differentiate \eqref{eq:overlap}.  Since the connection matrix is constant, we get
\[
B_{\ell'}(\tau)\Phi_{\ell'}(\tau)=B_\ell(\tau)\Phi_\ell(\tau)C_{\ell,\ell'}=B_\ell(\tau)\Phi_{\ell'}(\tau).
\]
Because $\Phi_{\ell'}$ is invertible, $B_{\ell'}=B_\ell$ on the overlap.  Thus the local coefficient matrices define a single meromorphic matrix $B$ on the union of the $U_\ell$.  The final assertion is the standard uniqueness statement for fundamental matrices: if $\Phi$ and $\Psi$ solve the same system on a simply connected domain, then $\Phi^{-1}\Psi$ has derivative zero, and hence is constant.
\end{proof}

\subsection[Proof of the modularity criterion]{Proof of Theorem~\ref{thm:MainResult}}

We assume that
 $\mathbf f=(f_1,\ldots,f_r)^{\mathsf T}$ is \emph{good} with exponents
$(\alpha_1,\dots,\alpha_r)$ for
$\mathcal O=\{(\gamma_\ell,w_\ell)\}_{\ell=1}^m$.

\medskip
\noindent
\emph{Finitely many orbits.}
For each $(\gamma_\ell,w_\ell)\in\mathcal O$, condition (i) gives
 $g_{\ell,1},\dots,g_{\ell,r}\in \SL_2(\mathbb Z)$ such that
\[
\Phi_\ell(\tau):=\big[\ \mathbf f(\gamma_\ell^{-1}g_{\ell,1}^{-1}\tau)
\ \cdots\ \mathbf f(\gamma_\ell^{-1}g_{\ell,r}^{-1}\tau)\ \big]
\]
is invertible on a punctured neighborhood of $q_\ell=0$ (\(q_\ell=e^{2\pi i\tau_\ell/w_\ell}\)).

\medskip
\noindent
\emph{Local first-order systems.}
By condition (ii), the matrices
\[
A_\ell(q_\ell):=(\theta_\ell\Phi_\ell)\,\Phi_\ell^{-1}
\qquad(\theta_\ell:=q_\ell\tfrac{d}{dq_\ell})
\]
are meromorphic on $0<|q_\ell|<1$, holomorphic at $q_\ell=0$, and satisfy
\(A_\ell(0)\sim\mathrm{diag}(\alpha_1,\ldots,\alpha_r)\).
In particular, each column of \(\Phi_\ell\) has the  form
\(q_\ell^{\Lambda}\cdot (\text{holomorphic invertible})\) with \(\Lambda=\mathrm{diag}(\alpha_j)\).
Furthermore, since $A_{\ell}(q_{\ell})$ is meromorphic on $0<|q_{\ell}| <1$ and the local coordinate $q_{\ell}=e^{2\pi i\tau_{\ell}/w_{\ell}}$ maps the chosen cusp neighborhood to a punctured disc, the induced system is defined on that neighborhood away from the isolated poles of the coefficient matrix.

\medskip

\noindent
{\it Orbit compatibility and finite coverage}. 
By condition~\textup{(iii)}, translating
$(\gamma_\ell,w_\ell)$ by $T^n$ yields a system that is gauge-equivalent (by a constant
matrix in $\mathrm{GL}_r(\mathbb C)$) to the original one, with the same limiting exponent
matrix. Since $\mathcal O$ is finite, all cuspidal patches that arise when we follow any path obtained from a word in $S$ and $T$ are modeled, up to constant gauge, on one of finitely many base systems indexed by $\mathcal{O}.$ Along any such path from
$\tau$ to $\gamma^{-1}\tau$ we therefore encounter only finitely many patches, and on each of them the first-order systems extend analytically.
 
\medskip
\noindent\emph{Global first-order system and connection matrices.}
For each $(\gamma_\ell,w_\ell)\in \mathcal O$, condition~(ii) gives a local
$q_\ell$-differential system
\[
\theta_\ell \Phi_\ell(\tau) = A_\ell(q_\ell)\,\Phi_\ell(\tau),
\]
whose columns are built from the $f(\gamma^{-1}\tau)$ and have
prescribed exponent matrix $\Lambda = \diag(\alpha_1,\dots,\alpha_r)$ at the
cusp corresponding to $\gamma_\ell$.  Passing from the coordinate $q_\ell$ to
$\tau$ on a simply connected neighborhood $U_\ell$ of that cusp, this induces a
first-order system of the form
\[
\frac{1}{2\pi i}\,\frac{d}{d\tau}\,\Phi_\ell(\tau)
  = B_\ell(\tau)\,\Phi_\ell(\tau),
\]
with $B_\ell(\tau)$ meromorphic on $U_\ell$.  By construction, on any overlap
$U_\ell\cap U_{\ell'}$ the columns of $\Phi_\ell$ and $\Phi_{\ell'}$ are analytic
continuations of the same $r$-tuples of translates of $f$, so the hypotheses of
Lemma~\ref{lem:gluing} are satisfied.  Therefore, there is a single globally
defined meromorphic matrix $B(\tau)$ on $\mathbb{H}$ such that
\begin{equation}\label{eq:global-system}
\frac{1}{2\pi i}\,\frac{d}{d\tau}\,\Phi_\ell(\tau)
  = B(\tau)\,\Phi_\ell(\tau)
\end{equation}
whenever $\Phi_\ell$ is defined, and on each overlap $U_\ell\cap U_{\ell'}$ we
have
\[
\Phi_{\ell'}(\tau) = \Phi_\ell(\tau)\,C_{\ell,\ell'},
\]
where $C_{\ell,\ell'}\in\GL_r(\C)$ is a constant connection matrix.

\medskip
\noindent
{\it Monodromy and transformation laws.}  Fix a base index $\ell_0$ such that $(\gamma_{\ell_0},w_{\ell_0})$
corresponds to the cusp at~$\infty$, and write $\Phi_\infty := \Phi_{\ell_0}$.  For any word
$\gamma$ in the generators $S$ and $T$, analytically continue the columns of $\Phi_\infty$
along a path from $\tau$ to $\gamma^{-1}\tau$ that winds only around cusps.  On each segment
of this path, the continuation is governed by the same global system
\eqref{eq:global-system}, and each time we cross from the domain of one $\Phi_\ell$ to that
of another $\Phi_{\ell'}$ we right-multiply by the corresponding constant connection matrix
$C_{\ell,\ell'}$.  The outcome is that
\[
\mathbf f(\gamma^{-1}\tau) \;=\; \rho(\gamma)\,\mathbf f(\tau),
\]
where $\rho(\gamma)$ is the product of the connection matrices encountered along the lift of
$\gamma$.  By uniqueness of analytic continuation, this $\rho(\gamma)$ depends only on
$\gamma$ and not on the chosen path or decomposition into patches, and the rule
$\gamma\mapsto \rho(\gamma)$ is multiplicative:
\[
\rho(\gamma_1\gamma_2)=\rho(\gamma_1)\rho(\gamma_2).
\]
In other words, $\rho$ is a representation of $\langle S,T\rangle$ into $\mathrm{GL}_r(\mathbb C)$,  and we have a vector-valued modular function.
Finally, taking $\gamma=T$ and using the expansion at the cusp at~$\infty$, we have
\[
\rho(T) = \lim_{q_\infty\to 0} e^{2\pi i\Lambda}
= \mathrm{diag}(e^{2\pi i\alpha_1},\ldots,e^{2\pi i\alpha_r}),
\]
as claimed in (A).  By definition, the transformation law
\[
\mathbf f(\gamma\tau)\;=\;\rho(\gamma)\,\mathbf f(\tau)\qquad(\gamma\in\SL_2(\mathbb Z))
\]
is precisely the condition that $\mathbf f$ be a vector-valued modular function with multiplier
$\rho$.  If this representation has finite image, then the components are ordinary weakly holomorphic modular functions on the finite-index subgroup $\ker\rho$; if $\rho$ factors through $\SL_2(\mathbb Z/N\mathbb Z)$, this is the corresponding congruence-level statement.

\medskip
\noindent
\emph{Weak holomorphy.}
Each \(f_j\) is holomorphic on \(\mathbb H\) by hypothesis.  At every cusp represented in the orbit datum, the corresponding fundamental matrix has the form \(q_\ell^{\Lambda}\) times a holomorphically invertible matrix, and hence each component has a Laurent expansion in the local parameter \(q_\ell\).  Therefore, \(\mathbf f\) is \emph{weakly holomorphic}.

\smallskip
\noindent
This completes the proof of
\emph{(B)$\ \Longrightarrow$ (A)}.

\medskip
\noindent
\emph{(A)$\ \Longrightarrow$ (B)}. 
If \(\mathbf f\) is vector-valued modular, then at each cusp \(c\) one
chooses \(r\) translates so that the columns are independent (exactly as we did in the case of the proof of Theorem~\ref{thm:RR-vvmf}), defines \(\Phi_c\) from those columns, and sets \(A_c=(\theta_c\Phi_c)\Phi_c^{-1}\).
This produces the matrices \(A_c\) and the finite orbit datum \(\mathcal O\) satisfying
the definition of ``goodness.''  The proof of Theorem~\ref{thm:converse} offers the explicit details.
$\qquad \square$
\medskip

\subsection[Sketch of the converse theorem]{Sketch of the Proof of Theorem~\ref{thm:converse}}
Assume \(\mathbf f=(f_1,\dots,f_r)^{\mathsf T}\) is a weakly holomorphic vector-valued modular function for a fixed congruence subgroup \(\Gamma\), with
\[
\mathbf f(\gamma\tau)=\rho(\gamma)\mathbf f(\tau)
\qquad\text{and}\qquad
\rho(T)=\mathrm{diag}(e^{2\pi i\alpha_1},\dots,e^{2\pi i\alpha_r}).
\]

\medskip
\noindent
\emph{Cusps and local parameters.}
Choose cusp representatives \(c=\gamma_c^{-1}\infty\in\Gamma\backslash\PP^1(\Q)\) with widths \(w_c\).
Then write
\[
\tau_c=\frac{a_c\tau+b_c}{c_c\tau+d_c},\qquad
q_c:=e^{2\pi i\,\tau_c/w_c},\qquad
\theta_c:=q_c\frac{d}{dq_c},\qquad
\mathbf f_c(\tau):=\mathbf f(\gamma_c^{-1}\tau).
\]
Weak holomorphy gives convergent Laurent expansions in the local parameter after the cusp monodromy is put in diagonal form.  Equivalently, for the chosen fundamental matrix of translates we may write
\[
   \Phi_c(\tau)=H_c(q_c)\,q_c^\Lambda C_c,
   \qquad \Lambda=\diag(\alpha_1,\dots,\alpha_r),
\]
where $H_c(q_c)$ is holomorphic and invertible at $q_c=0$, and $C_c$ is a constant invertible matrix.  This is the usual Frobenius form for a regular singular first-order system at a cusp; it follows from the finite-order diagonal $T$-monodromy and the weak holomorphy assumption.  In concrete examples, including the Rogers--Ramanujan pair, the columns are obtained by using suitable translates such as powers of $T$.

\medskip
\noindent
\emph{Fundamental columns.}
For each cusp \(c\), choose \(g_{c,1},\dots,g_{c,r}\in \SL_2(\Z)\) so that
\[
\Phi_c(\tau):=\big[\ \mathbf f(\gamma_c^{-1}g_{c,1}^{-1}\tau)\ \cdots\
                   \mathbf f(\gamma_c^{-1}g_{c,r}^{-1}\tau)\ \big]
\]
has the Frobenius form above.  Since $H_c(0)$ and $C_c$ are invertible, this matrix is invertible for $q_c$ sufficiently small.

\smallskip
\noindent
\emph{Local $q$-differential system.}
We let
\[
A_c(q_c):=(\theta_c\Phi_c(\tau))\,\Phi_c(\tau)^{-1}.
\]
As the entries of \(\Phi_c\) are meromorphic in \(q_c\) and \(\det\Phi_c\not\equiv0\), we find that
\(A_c\) is meromorphic on \(0<|q_c|<1\).
From the Frobenius form \(\Phi_c=H_c(q_c)q_c^{\Lambda}C_c\), we have
\[
A_c(q_c)=(\theta_c H_c)\,H_c^{-1}+H_c(q_c)\Lambda H_c(q_c)^{-1}.
\]
This extends holomorphically to \(q_c=0\), and since \(\theta_c H_c\) has zero constant term,
\[
A_c(0)=H_c(0)\Lambda H_c(0)^{-1}\sim\mathrm{diag}(\alpha_1,\dots,\alpha_r).
\]
Therefore, \(\theta_c\Phi_c=A_c(q_c)\Phi_c\) gives the desired  $q$-differential system.

\medskip
\noindent
\emph{Orbit compatibility and finiteness.}
Replacing \(\gamma_c\) by \(\gamma_c T^n\) multiplies each column of \(\Phi_c\) by the constant matrix
\(\rho(T)^{-n}\). Therefore, the resulting system is gauge-equivalent with the same limiting exponents.
Since there are finitely many cusps and we chose \(r\) columns at each, the orbit datum
\[
\mathcal O=\bigcup_{\,c=\gamma_c^{-1}\infty}\{(\gamma_c g_{c,1},w_c),\dots,(\gamma_c g_{c,r},w_c)\}
\]
is finite and satisfies the $T$-stability in Definition~\ref{def:good} (iii).

\smallskip
\noindent
To summarize,
we have produced, for each cusp \(c\), an invertible \(\Phi_c\) (i.e. condition (i)),
a logarithmic $q$-differential system \(\theta_c\Phi_c=A_c(q_c)\Phi_c\) with \(A_c\) holomorphic at \(q_c=0\) and
\(A_c(0)\sim\mathrm{diag}(\alpha_j)\) (i.e. condition (ii)), and the orbit-compatibility/finite coverage
(i.e. condition (iii)). Therefore, \(\mathbf f\) is \emph{good} for the explicit finite set \(\mathcal O\). $\qquad \square$

\section[Andrews--Gordon illustration]{The Andrews--Gordon family}\label{sec:AG}
Here we indicate how the Rogers--Ramanujan argument extends to the normalized Andrews--Gordon vectors $\mathbf f_k$.  This section is intended as a guide to the method rather than as a replacement for the detailed finite polynomial identities in the Andrews--Gordon literature.

Schur's recurrences are generalized by the Andrews--Gordon polynomial method.  The required finite polynomial identities are supplied by Andrews, Gordon, and Warnaar \cite{Andrews1974,Gordon1961,Warnaar}.  They produce finite truncations whose limits are the series $F_{k,i}$ and whose product limits are the products in \eqref{AG}.  As in the Rogers--Ramanujan case, the logarithmic derivatives of these finite identities stabilize coefficientwise and yield first-order systems for the normalized vector $\mathbf f_k$.

To describe the limiting coefficient functions, set $M=2k+1$ and define the weighted residue Lambert series
\[
L_{r,M}(q):=\sum_{\substack{n\ge1\\ n\equiv r\pmod M}}\frac{nq^n}{1-q^n}
\qquad (r\in\mathbb Z/M\mathbb Z).
\]
Direct logarithmic differentiation of the Andrews--Gordon product gives
\[
\theta\log F_{k,i}(q)
=\sum_{n\ge1}\frac{nq^n}{1-q^n}
 -L_{0,M}(q)-L_{i,M}(q)-L_{M-i,M}(q).
\]
Equivalently, after the normalization by $q^{\alpha_{k,i}}$, the diagonal entries of the resulting coefficient matrix are
\[
\alpha_{k,i}+\sum_{n\ge1}\frac{nq^n}{1-q^n}
 -L_{0,M}(q)-L_{i,M}(q)-L_{M-i,M}(q).
\]
The finite Andrews--Gordon recurrences give the corresponding finite systems before passage to the limit; the displayed formula records the limiting coefficient functions needed for Definition~\ref{def:good}.

The same cuspwise argument as in Section~\ref{RRProof} then applies.  The $T$-matrix is read off from the exponents $\alpha_{k,i}$.  The $S$-matrix is the standard sine matrix for the $(2,2k+1)$ minimal-model character vector, equivalently the connection matrix for the Andrews--Gordon system,
\[
\big(\rho_k(S)\big)_{ij}
=\frac{2}{\sqrt{2k+1}}\,\sin\!\left(\frac{2\pi i j}{2k+1}\right).
\]
This is the standard transformation law for the $(2,2k+1)$ minimal-model character vector, up to the diagonal gauge and ordering conventions mentioned above.  The same modularity is also a consequence of Zhu's theorem on vertex-operator-algebra characters \cite{Zhu1996}; here the point is that the finite Andrews--Gordon polynomial identities provide the $q$-series differential systems used by Theorem~\ref{thm:MainResult}.

\end{document}